\documentclass[a4paper,11pt]{amsart}
\usepackage[colorlinks, linkcolor=blue,anchorcolor=Periwinkle,
    citecolor=blue,urlcolor=Emerald]{hyperref}
\usepackage[all]{xy}
\SelectTips{cm}{}

\usepackage{graphicx}
\usepackage{psfrag}

\usepackage{mathtools}
\usepackage{tikz}
\usepackage{tikz-cd}
\usepackage{extarrows}

\usepackage{amssymb}
\usetikzlibrary{decorations.pathreplacing}
\usetikzlibrary{matrix,arrows}

\newcommand{\ie}{{i.e.}\ }
\newcommand{\cf}{{cf.}\ }

\newcommand{\ko}{\: , \;}
\newcommand{\ul}[1]{\underline{#1}}
\newcommand{\ol}[1]{\overline{#1}}

\newtheorem{theorem}{Theorem}
\newtheorem{classification-theorem}[subsection]{Classification Theorem}
\newtheorem{decomposition-theorem}[subsection]{Decomposition Theorem}
\newtheorem{proposition-definition}[subsection]{Proposition-Definition}
\newtheorem{periodicity-conjecture}[subsection]{Periodicity Conjecture}
\newtheorem{lemma}[theorem]{Lemma}

\newtheorem{example}[theorem]{Example}
\newtheorem{remark}[theorem]{Remark}

\newcommand{\reminder}[1]{}

\renewcommand{\mod}{\mathrm{mod}\,}
\newcommand{\sg}{\mathrm{sg}\,}

\newcommand{\Mod}{\mathrm{Mod}\,}
\newcommand{\proj}{\mathrm{proj}\,}

\newcommand{\per}{\mathrm{per}\,}
\newcommand{\pvd}{\mathrm{pvd}\,}
\newcommand{\add}{\mathrm{add}\,}

\newcommand{\Tr}{\mathrm{Tr}}

\newcommand{\Si}{\Sigma}
\newcommand{\Ga}{\Gamma}

\newcommand{\C}{\mathbb{C}}
\newcommand{\R}{\mathbb{R}}

\newcommand{\iso}{\xrightarrow{_\sim}}

\newcommand{\id}{\mathbf{1}}

\newcommand{\Hom}{\mathrm{Hom}}

\newcommand{\End}{\mathrm{End}}

\newcommand{\ten}{\otimes}

\newcommand{\ca}{{\mathcal A}}
\newcommand{\cb}{{\mathcal B}}
\newcommand{\cc}{{\mathcal C}}
\newcommand{\cd}{{\mathcal D}}
\newcommand{\ce}{{\mathcal E}}
\newcommand{\cF}{{\mathcal F}}

\newcommand{\ch}{{\mathcal H}}

\newcommand{\cm}{{\mathcal M}}
\newcommand{\cn}{{\mathcal N}}

\newcommand{\cp}{{\mathcal P}}

\newcommand{\ct}{{\mathcal T}}

\begin{document}



\title{Calabi--Yau structures on Drinfeld quotients and Amiot's conjecture}
\author{Bernhard Keller}
\address{Universit\'e Paris Cit\'e and Sorbonne Université, CNRS, IMJ-PRG, F-75013 Paris, France}
\email{bernhard.keller@imj-prg.fr}
\urladdr{https://webusers.imj-prg.fr/~bernhard.keller/}

\author{Junyang Liu}
\address{Universit\'e Paris Cit\'e and Sorbonne Université, CNRS, IMJ-PRG, F-75013 Paris, France 
and Yau Mathematical Sciences Center, Tsinghua University, Beijing 100084, China}
\email{liuj@imj-prg.fr}
\email{liujunya18@mails.tsinghua.edu.cn}
\urladdr{https://webusers.imj-prg.fr/~junyang.liu}

\begin{abstract}
In 2009, Claire Amiot gave a construction of Calabi--Yau structures
on Verdier quotients. We sketch how to lift it to the dg setting. 
We use this construction as an important step in an outline of 
the proof of her conjecture on the structure of $2$-Calabi--Yau 
triangulated categories with a cluster-tilting object.
\end{abstract}





\maketitle

In 2009, Amiot \cite{Amiot09} constructed Calabi--Yau structures on certain
Verdier quotients. Our first aim is to lift her construction to the level of dg (=differential
graded) enhancements. Our second aim is to sketch how this construction
allows one to prove Amiot's conjecture \cite{Amiot11} on the structure of
$2$-Calabi--Yau triangulated categories with a cluster-tilting object.

\section{Calabi--Yau structures on Drinfeld quotients}
Let $k$ be a field and $\cn$ a $k$-linear $\Hom$-finite triangulated category.
Let $d$ be an integer and suppose that $\cn$ is $d$-Calabi--Yau, \ie endowed with
bifunctorial isomorphisms
\[
\Hom(N_1, N_2) \xlongrightarrow{_\sim} D\Hom(N_2, \Si^d N_1) \: ,
\]
where $N_1$, $N_2$ lie in $\cn$ and $D$ denotes the dual over $k$. For example, let $X$ be a smooth
projective variety of dimension $d$ over an algebraically closed ground field
and $\cn$ the bounded
derived category of coherent sheaves on $X$. Then, by one possible definition, the
variety $X$ is $d$-Calabi--Yau if and only if its canonical sheaf $\omega_X$ is trivialisable
and that happens if and only if $\cn$ is $d$-Calabi--Yau as a triangulated category. As a second
example, suppose that $A$ is a finite-dimensional basic algebra and $\cn$ the bounded
homotopy category of finitely generated projective $A$-modules. Then $\cn$ is $0$-Calabi--Yau
if and only if $A$ is symmetric, \ie isomorphic to its dual $DA$ as a bimodule. The examples we are
most interested in are associated with quivers with potential $(Q,W)$ 
which are \linebreak Jacobi-finite, \ie the associated (complete)
Jacobian algebra is finite-dimensional, \cf \cite{DerksenWeymanZelevinsky08} for the terminology. 
Let $\Ga$ be the Ginzburg dg algebra \cite{Ginzburg06},
\cf also \cite{KellerYang11}, associated with such a quiver with potential. Let $\cd(\Ga)$ be
its unbounded derived category, $\per \Ga$ its perfect derived category (the full subcategory
of the compact objects in $\cd(\Ga)$) and $\pvd \Ga$ its {\em perfectly valued derived category},
\ie the full subcategory of objects $M$ of $\cd(\Ga)$ such that the underlying complex
$M|_k$ is perfect over $k$. By Amiot's definition \cite{Amiot09} of the cluster category
$\cc_{Q,W}$, we have an exact sequence
\[
\begin{tikzcd}
0 \ar[r] & \pvd \Ga \ar[r] & \per \Ga \ar[r] & \cc_{Q,W} \ar[r] & 0
\end{tikzcd}
\]
of triangulated categories. Here the perfectly valued derived category is $3$-Calabi--Yau, the perfect derived
category is not Calabi--Yau and the cluster category is $2$-Calabi--Yau, \cf \cite{Amiot09}
and the references given there.

The last example can be generalised as follows: suppose that the category $\cn$ is contained as a thick subcategory in a $k$-linear triangulated category $\ct$. As shown by Amiot \cite{Amiot09}, in the exact sequence
\[
\begin{tikzcd}
0 \ar[r] & \cn \ar[r] & \ct \ar[r] & \ct/\cn \ar[r] & 0 \: ,
\end{tikzcd}
\]
the Verdier quotient $\ct/\cn$ is often $(d-1)$-Calabi--Yau. Let us recall the
construction of the Calabi--Yau structure in a slightly more general setting:
let us assume that we are given a {\em bilinear form of degree $-d$ on
$\cn$}, \ie a family of bilinear maps
\[
\cn(N_1, N_2) \times \cn(N_2, \Si^d N_1) \longrightarrow k
\]
such that the associated maps
\begin{equation} \label{eq:bilin-form}
\cn(N_1, N_2) \longrightarrow D\cn(N_2, \Si^d N_1)
\end{equation}
define a morphism of bifunctors on $\cn$. 
Such a bilinear form is determined by the {\em pretrace forms}
\[
t_N\colon \cn(N, \Si^d N) \longrightarrow k
\]
obtained as the images of the identities $\id_N$ under the maps~(\ref{eq:bilin-form}). 
One defines a {\em $d$-Calabi--Yau structure on $\cn$} to be a bilinear form of degree $-d$ on 
$\cn$ such that the map~(\ref{eq:bilin-form}) is invertible and which satisfies a graded symmetry property, \cf Proposition~A.5.2 of \cite{Bocklandt08}.
\begin{theorem}[Amiot \cite{Amiot09}] The category $\ct/\cn$ admits a unique
bilinear form of degree $1-d$ whose pretrace form $t_X^{\ct/\cn}$ for an
object $X$ is given by 
\[
t_X^{\ct/\cn}(``f \circ s^{-1}") = t_N^{\cn}((\Si^{d-1} b) \circ f \circ a) \: ,
\]
where $``f\circ s^{-1}"$ is a fraction given by morphisms $f$ and $s$ (with cone $\Si N$ in $\cn$)
as in the following diagram
\[
\begin{tikzcd}
N \ar[r,"a"] & X' \ar[r,"s"] \ar[d, "f"] & X \ar[r,"b"] & \Si N \\
                 & \Si^{d-1} X \ar[r, "\Si^{d-1} b"] & \Si^d N\mathrlap{\: .} & 
\end{tikzcd}
\]
Moreover, if the bilinear form on $\cn$ is a $d$-Calabi--Yau structure and $\cn\subseteq \ct$
satisfies suitable non-degeneracy conditions, then the bilinear form on
$\ct/\cn$ thus constructed is a $(d-1)$-Calabi--Yau structure.
\end{theorem}
For example, if $A$ is a finite-dimensional symmetric algebra, we have the short exact
sequence
\[
\begin{tikzcd}
0 \ar[r] & \ch^b(\proj A) \ar[r] & \cd^b(\mod A) \ar[r] & \sg A \ar[r] & 0
\end{tikzcd}
\]
of triangulated categories and we find that the stable category $\ul{\mathrm{mod}}\, A\iso \sg A$ is \linebreak $(-1)$-Calabi--Yau,
which corresponds to the well-known fact that the Auslander--Reiten translation is given by the square $\Omega^2=\Si^{-2}$ of the syzygy functor.

Our aim is to lift this construction to the level of dg categories.
Let $\ca$ be a (small) dg category and $\cd(\ca)$ its derived category. Its objects are the
dg functors $M\colon \ca^{op} \to \cc_{dg}(\Mod k)$ from $\ca^{op}$ to the dg category
of complexes over $k$. Suppose that $\ca$ is {\em pretriangulated}, \ie the
Yoneda functor
\[
H^0(\ca) \to \cd(\ca) \ko X \mapsto \ca(?,X)
\]
is an equivalence onto a full triangulated subcategory. Let $\cb\subseteq\ca$ be a pretriangulated full dg subcategory. By definition \cite{Drinfeld04}, the {\em Drinfeld quotient}
$\ca/_{Dr} \,\cb = \ca/\cb$ is obtained from $\ca$ by formally adjoining a contracting
homotopy
\[
h_N\colon N \to N \ko |h_N|=-1 \ko d(h_N)= \id_N
\]
for each object $N$ in $\cb$. Thus, by definition, the canonical functor $\ca \to \ca/\cb$ is
{\em strictly universal} among the dg functors $F\colon \ca \to \cF$ to a dg category where
all the objects $FN$, $N$ in $\cb$, are endowed with a contracting homotopy. 
As shown in \cite{Drinfeld04}, we have an induced short exact sequence
\[
\begin{tikzcd}
0 \ar[r] & H^0(\cb) \ar[r] & H^0(\ca) \ar[r] & H^0(\ca/\cb) \ar[r] & 0
\end{tikzcd}
\]
of triangulated categories, \cf~also \cite{Keller99}. It follows that the sequence
\[
\begin{tikzcd}
0 \ar[r] & \cb \ar[r] & \ca \ar[r] & \ca/\cb \ar[r] & 0
\end{tikzcd}
\]
of dg categories is {\em homotopy exact}, \ie the sequence
\[
\begin{tikzcd}
0 \ar[r] &\cd(\cb) \ar[r] & \cd(\ca) \ar[r] & \cd(\ca/\cb) \ar[r] & 0
\end{tikzcd}
\]
of triangulated categories is exact. Therefore, by the main result of \cite{Keller98}, we have long exact
sequences
\[
\begin{tikzcd}
\cdots \ar[r] & HC_*(\cb) \ar[r] & HC_*(\ca) \ar[r] & HC_*(\ca/\cb) \ar[r] & HC_{*-1}(\cb) \ar[r] & \cdots
\end{tikzcd}
\]
in Hochschild and cyclic homology. Now suppose that $\cb$ is an {\em $H^*$-finite} dg category, \ie the spaces 
$H^p(\cb)(X,Y)$ are
finite-dimensional for all $X$, $Y$ in $\cb$ and all integers $p$.
Let $d$ be an integer. Following Kontsevich, a (right) {\em $d$-Calabi--Yau structure} on
$\cb$ is a class $c$ in $DHC_{-d}(\cb)$ which is {\em non-degenerate}, \ie its image under
the composition
\[
\begin{tikzcd}
DHC_{-d}(\cb) \arrow{r} & DHH_{-d}(\cb) \arrow[no head]{r}{\sim} & \Hom_{\cd(\cb^e)}(\cb, \Si^{-d}D\cb)
\end{tikzcd}
\]
of canonical maps is invertible. Here, by abuse of notation, we denote by $\cb$ the $\cb$-bimodule
$(X,Y) \mapsto \cb(X,Y)$ and by $D\cb$ the $\cb$-bimodule $(X,Y) \mapsto D\cb(Y,X)$.
If $\cb$ is pretriangulated and carries a $d$-Calabi--Yau structure, then $H^0(\cb)$ becomes
$d$-Calabi--Yau as a triangulated category since we have
\[
\begin{tikzcd}
H^0(\cb)(X,Y) \arrow[no head]{r}{\sim} & H^0(\Si^{-d} D\cb(Y,X)) \arrow[no head]{r}{\sim} & D H^0(\cb)(Y, \Si^d X) \: .
\end{tikzcd}
\]
More generally, we obtain a canonical map
\[
\begin{tikzcd}
DHC_{-d}(\cb) \ar[r] & \{ \mbox{bilinear forms of degree $-d$ on $H^0(\cb)$}\} \: ,
\end{tikzcd}
\]
which induces a map
\[
\begin{tikzcd}
DHC_{-d}^{nd}(\cb) \ar[r] & \{ \mbox{$d$-Calabi--Yau structures on $H^0(\cb)$}\} \: ,
\end{tikzcd}
\]
where $DHC_{-d}^{nd}(\cb)$ denotes the space of non-degenerate classes in $DHC_{-d}(\cb)$.

Now let $\cb\subseteq \ca$ be a pretriangulated full dg subcategory of a pretriangulated dg category $\ca$.
Thus, we have a long exact sequence
\[
\begin{tikzcd}
\cdots \ar[r] & HC_{1-d}(\cb) \ar[r] & HC_{1-d} (\ca) \ar[r] & HC_{1-d}(\ca/\cb) \ar[r,"\delta"] & 
HC_{-d} (\cb) \ar[r] & \cdots\: .
\end{tikzcd}
\]
\begin{theorem}[Keller--Liu \cite{KellerLiu23}] \label{theorem:connecting morphism} The square
\[
\begin{tikzcd}
DHC_{-d}(\cb) \ar[d] \ar[r, "D\delta"] & DHC_{1-d}(\ca/\cb) \ar[d] \\
 \{ \mbox{bilin.~forms of degree $-d$ on $H^0(\cb)$}\} \ar[r] &
  \{ \mbox{bilin.~forms of degree $1-d$ on $H^0(\ca/\cb)$}\}\mathrlap{\: ,}
  \end{tikzcd}
\]
where the bottom horizontal map is given by Amiot's construction,
is commutative.
\end{theorem}
Let us sketch the proof: let $HH(\ca)$ be the Hochschild complex of $\ca$.
It is the sum total complex of the bicomplex
\[
\begin{tikzcd}
\coprod_{A_0} \ca(A_0, A_0) & \coprod_{A_0, A_1} \ca(A_0, A_1)\ten \ca(A_1,A_0) \ar[l] &
\cdots \: , \ar[l]
\end{tikzcd}
\]
where, for example, the leftmost
differential sends $f_0\ten f_1$ to
\[
(-1)^{|f_0|}f_0 \circ f_1 - (-1)^{|f_0| (|f_1|+1)} f_1 \circ f_0 \: .
\]
Then the sequence
\[
\begin{tikzcd}
0 \ar[r] &\cb \ar[r] & \ca \ar[r] & \ca/\cb \ar[r] & 0
\end{tikzcd}
\]
of dg categories with the given null-homotopy for the composition $\cb \to \ca/\cb$ yields a {\em homotopy
short exact sequence of complexes} (defined below)
\begin{equation} \label{HH-sequence}
\begin{tikzcd}
HH(\cb) \ar[rr, "h", bend left] \ar[r, "i", swap] & HH(\ca) \ar[r, "p", swap] & HH(\ca/\cb)\: ,
\end{tikzcd}
\end{equation}
where the map $h$ sends $f_0$ to
$h_{B_0} f_0$ ($B_0$ is the source and the target of $f_0$). The {\em homotopy snake lemma} (stated below) applied to the
sequence (\ref{HH-sequence}) then allows
us to make the connecting morphism 
\[
\delta \colon HH_{1-d}(\ca/\cb) \longrightarrow HH_{-d}(\cb)
\]
explicit and to check the commutativity claimed in the theorem.

We define a {\em homotopy short exact sequence of complexes} to be a diagram
of complexes of abelian groups
\[
\begin{tikzcd}
B \ar[r, "i", swap] \ar[rr, "h", bend left] & A \ar[r, "p", swap] & C
\end{tikzcd}
\]
such that $i$ and $p$ are morphisms of complexes and $h$ is a homogeneous
morphism of degree $-1$ such that we have $d(h) = p\circ i$ and that the graded
object $\Si^{-1} C \oplus A \oplus \Si B$ endowed with 
the differential
\[
\begin{bmatrix}
-d & i & -h \\
0 & d & -p \\
0 & 0 & -d
\end{bmatrix}
\]
is acyclic. Such a diagram yields a triangle
\[
\begin{tikzcd}
B \ar[r] & A \ar[r] & C \ar[r] & \Si B
\end{tikzcd}
\]
in the derived category of abelian groups and we can compute
the connecting morphism
\[
\delta \colon H^q(C) \to H^{q+1}(B)
\]
thanks to the following homotopy snake lemma, where $a\in A^q$, $b\in B^{q+1}$, $c\in C^q$.
\begin{lemma}
We have $\delta(\ol{c}) = -\ol{b}$ if and only if there is an element $a$ in $A$ such that we have $d(a) + i(b) = 0$ and $\ol{p(a) + h(b)} =\ol{c}$.
\end{lemma}
As mentioned above, we apply this lemma to the sequence~(\ref{HH-sequence}) to 
prove the commutativity stated in the Theorem~\ref{theorem:connecting morphism}.

\section{Amiot's conjecture}
Suppose that $k$ is a field of characteristic $0$.
Let $\cc$ be a $k$-linear triangulated category which is $\Hom$-finite and Karoubian.
Assume that $\cc$ is $2$-Calabi--Yau. Recall that an object $T$ of $\cc$ is
{\em cluster-tilting} if
\begin{itemize}
\item[a)] we have $\cc(T, \Si T)=0$ and
\item[b)] for each object $X$ of $\cc$, there is a triangle
\[
\begin{tikzcd}
T_1 \ar[r] &  T_0 \ar[r] & X \ar[r] & \Si T_1
\end{tikzcd}
\]
with $T_i$ in $\add T$.
\end{itemize}
Amiot showed in \cite{Amiot09} that if $(Q,W)$ is a Jacobi-finite quiver with
potential, then the cluster category $\cc_{Q,W}$ is $k$-linear, ($\Hom$-finite,) Karoubian, $2$-Calabi--Yau
and admits a canonical cluster-tilting object. In Question~2.20 of \cite{Amiot11}, she asked whether,
under suitable additional hypotheses, all examples are obtained in this way. The question
is discussed for example in \cite{KalckYang16}. The fact that the
answer should be positive is now known as `Amiot's conjecture', \cf for example Conjecture~1.1 of 
\cite{KalckYang20}. Classes of examples can be found in 
\cite{AmiotReitenTodorov11,
AmiotIyamaReitenTodorov12,
Garcia21}.

\subsection{The quiver case}
\begin{theorem} \label{thm:quiver-case} Assume that $k$ is algebraically closed and
\begin{itemize}
\item[a)] there is a Frobenius exact category $\ce$ enriched over the category of pseudo-compact vector spaces with the full subcategory $\cp$ of projective-injective objects and the stable category $\cc$ such that
\item[b)] the dg category $\cc_{dg}=\cd^b_{dg}(\ce)/\cd^b_{dg}(\cp)$ carries a right $2$-Calabi--Yau structure inducing the given $2$-Calabi--Yau structure on $\cc$ and
\item[c)] $\cc$ contains a basic cluster-tilting object $T$.
\end{itemize}
Then there exists a quiver with potential $(Q,W)$ and a triangle equivalence $\cc_{Q,W}\iso\cc$ taking the canonical cluster-tilting object of $\cc_{Q,W}$ to $T$. 
Consequently, the endomorphism algebra of $T$ is isomorphic to a Jacobian algebra.
\end{theorem}
This theorem is a consequence of the more general Theorem~\ref{thm:species-case} stated
below, which also holds for ground fields $k$ which are not algebraically closed. 

\subsection{The species case}
Suppose that $l$ is a finite-dimensional semisimple $k$-algebra and $V_c$ a pseudo-compact graded $l$-bimodule concentrated in degrees $[-1,0]$ of finite total dimension. We denote by $T_l V_c$ the {\em completed} tensor graded algebra of $V_c$ over $l$ and by $[T_l V_c,T_l V_c]$ the closure of the vector subspace generated by the commutators.
Let $\eta$ be a non-degenerate and anti-symmetric element in $V_c\ten_{l^e}V_c$ of degree $-1$ and $w$ an element in $T_l V_c/[T_l V_c,T_l V_c]$ of degree $0$. 
Following Van den Bergh \cite{VandenBergh15}, we define 
the {\em deformed dg preprojective algebra} $\Pi(l, V_c,\eta,w)$ associated with the data $(l, V_c,\eta,w)$ as the $l$-augmented pseudo-compact dg algebra $T_l(V_c\oplus zl)$, where $z$ is an $l$-central element of degree $-2$
and the differential is determined by
\[
d(z)=\sigma'\eta\sigma''\quad \mbox{and}\quad d(f)=\{w,f\}_{\omega_\eta}\quad \mbox{for all }f \in T_l V_c.
\]
Here the symbol $\sigma=\sigma'\ten \sigma''$ denotes the Casimir element in $l\ten_k l$ associated with a fixed trace
form on $l$ and $\{?,-\}_{\omega_\eta}$ is the necklace bracket associated with the bisymplectic form 
$\omega_\eta$ defined by $\eta$, \cf section 10.1 of \cite{VandenBergh15}.

Let us make the assumptions a), b) and c) of the preceding theorem but not assume that $k$ is algebraically
closed or that $T$ is basic. 
Let $l$ be the largest semisimple quotient of $\End(T)$. Since $k$ is of characteristic $0$, the
semisimple $k$-algebra $l$ is separable. Therefore, by the Wedderburn--Mal'cev theorem, we can
choose an augmentation for $\End(T)$ over $l$. 
\begin{theorem} \label{thm:species-case}
There exists a deformed dg preprojective algebra $\Pi=\Pi(l, V_c,\eta,w)$ and a triangle equivalence 
$\cc_\Pi\iso\cc$ taking the canonical cluster-tilting object of $\cc_\Pi$ to $T$.
Consequently, the endomorphism algebra of $T$ is isomorphic to $H^0(\Pi)$.
\end{theorem}
We can deduce Theorem~\ref{thm:quiver-case} from Theorem~\ref{thm:species-case}:
indeed, when the ground field $k$ is algebraically closed and $T$ is basic, the algebra $l$ is isomorphic
to a finite product of copies of the ground field and $l$-bimodules can be described using
vector spaces with bases given by arrows. The `deformation' contained in the term
`deformed dg preprojective algebra' is then encoded in the potential.

Let us sketch the proof of Theorem~\ref{thm:species-case}:
we define $\Ga$ to be the endomorphism algebra of the image of $T$ in the dg quotient $\cc^b_{dg}(\cm)/\cc^b_{dg}(\cp)$, where $\cm\subseteq \ce$ denotes the closure under finite direct sums and direct summands of $T$ and the projective-injective objects in $\ce$. It is an $l$-augmented pseudo-compact dg algebra.
One can show that $\cc_{dg}$ fits into a short exact sequence
\[
\begin{tikzcd}
0 \ar[r] & \pvd\!_{dg}\Ga \ar[r] & \per\!_{dg}\Ga \ar[r] & \cc_{dg} \ar[r] & 0
\end{tikzcd}
\]
of dg categories. Therefore, it suffices to show that $\Ga$ is quasi-isomorphic to a
deformed dg preprojective dg algebra $\Pi$.
To prove this using the main result of \cite{VandenBergh15}, it suffices to construct a
non-degenerate class in negative cyclic homology $HN_3(\Ga)$. Since $\pvd\!_{dg}\Ga$
is quasi-equivalent to $\per\!_{dg}\Ga^!$, where $\Ga^!$ is the Koszul dual of $\Ga$,
the above exact sequence of dg categories
yields \cite{Keller98}  a long exact sequence in cyclic homology with invertible connecting morphism
\[
\delta \colon HC_{-2}(\cc_{dg}) \longrightarrow HC_{-3}(\Ga^!) \: .
\]
Moreover, by Corollary~D.2 of \cite{VandenBergh15}, we have a canonical $k$-linear bijection
\[
HN_3(\Ga) \xlongrightarrow{_\sim} DHC_{-3}(\Ga^!) \: .
\]
Thus, we get a composed bijection
\[
\begin{tikzcd}
HN_3(\Ga) \ar[r,"_\sim"] &  DHC_{-3}(\Ga^!) \ar[r,"D\delta"] & DHC_{-2}(\cc_{dg}) \: .
\end{tikzcd}
\]
We let $\alpha$ be the given non-degenerate class in $DHC_{-2}(\cc_{dg})$ and $\gamma$ its
preimage in $HN_3(\Ga)$. The subtle point is to prove that $\gamma$ is non-degenerate. 
The key to proving this is Theorem~\ref{theorem:connecting morphism}. More precisely,
we proceed as follows: it is not hard to see that the non-degeneracy of $\gamma$
follows if its image $\beta$ in $DHC_{-3}(\Ga^!)$ is non-degenerate.
For this, we show that $D\delta$ detects non-degeneracy. To prove this point, we use
Theorem~\ref{theorem:connecting morphism} to reduce the claim to a statement
about graded symmetric bilinear forms on the categories $\cc=H^0(\cc_{dg})$ and
$\pvd \Ga=H^0(\pvd\!_{dg}\Ga)$. We prove this statement by showing that the
construction in the proof of part~c) of Proposition~4 of \cite{KellerReiten07} is `locally inverse' to Amiot's construction.

\begin{example} This example serves to illustrate the construction of $\Pi(l, V_c,\eta,w)$ in the case where $k$ is not algebraically closed. Starting from the $\R$-species
\tikzset{
    lablp/.style={anchor=south, rotate=39, inner sep=.5mm},
    lablm/.style={anchor=south, rotate=-39, inner sep=.5mm},
}
\[
\begin{tikzcd}
	&\R\arrow[dr,"\R c\oplus\R ic" lablm]&\\
	\R\arrow[ur,"\R a" lablp]&&\C\arrow[ll,"\R b\oplus\R bi"]
\end{tikzcd}
\]
with the potential $W=abc$ we obtain the graded double $\R$-species
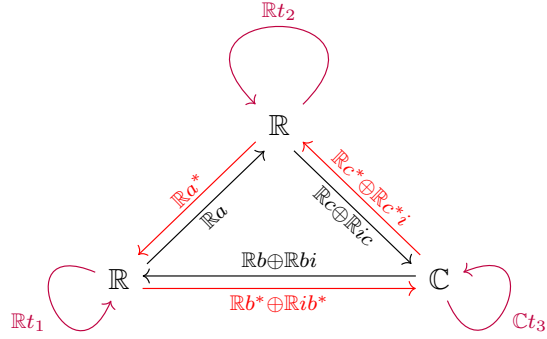
\tikzset{
    lablsouthp/.style={anchor=south, rotate=44, inner sep=.5mm},
    lablnorthp/.style={anchor=north, rotate=44, inner sep=.5mm},
    lablsouthm/.style={anchor=south, rotate=-44, inner sep=.5mm},
    lablnorthm/.style={anchor=north, rotate=-44, inner sep=.5mm},
}
\[
\begin{tikzcd}[row sep=1.5cm, column sep=1.5cm]
	&\R\arrow[dr,shift right=0.5ex,"\R c\oplus \R ic" lablnorthm,swap]
	\arrow[dl,shift right=0.5ex,"\R a^*" lablsouthp,swap,red]
	\arrow[out=45,in=135,loop,purple,"{\R t_{2}}",swap]&\\
	\R\arrow[ur,shift right=0.5ex,"\R a" lablnorthp,swap]
	\arrow[rr,shift right=0.5ex,"\R b^*\oplus \R ib^*",swap,red]
	\arrow[out=157,in=247,loop,purple,"{\R t_{1}}",swap]&&
	\C\arrow[ll,shift right=0.5ex,"\R b\oplus\R bi",swap]
	\arrow[ul,shift right=0.5ex,"\R c^*\oplus\R c^*i" lablsouthm,swap,red]
	\arrow[out=293,in=23,loop,purple,"{\C t_{3}}",swap]
\end{tikzcd}
\]
with the elements $a^*$, $b^*$, $c^*$ of degree $-1$ and $t_1$, $t_2$, $t_3$ of degree $-2$. It gives rise to the $l$-augmented pseudo-compact dg $\R$-algebra 
$\Pi=T_l(V_c\oplus\R t_1\oplus\R t_2\oplus\C t_3)$, where $l$ is the semisimple $\R$-algebra $\R\times\R\times\C$ and $V_c$ is the $l$-bimodule with the $\R$-basis
\[
a\ko b\ko bi\ko c\ko ic\ko a^*\ko b^*\ko ib^*\ko c^*\ko c^*i\: .
\]
In the notations
above, we have $z=t_1+t_2+t_3$ and $\sigma$ is the Casimir element associated with the trace form $\Tr(\lambda_1,\lambda_2,\lambda_3)=\lambda_1+\lambda_2+\mathrm{Re}(\lambda_3)$.

For a suitable choice of non-degenerate and graded anti-symmetric element $\eta$ in $V_c\ten_{l^e}V_c$, we obtain
the differential determined by
\[
d(t_1)=bb^*+bib^*-a^*a \ko
d(t_2)=aa^*-c^*c-c^*ic \ko
d(t_3)=(1+i)(cc^*-b^*b)+(1-i)(cc^*-b^*b)i
\]
and
\[
d(a^*) = -bc \ko
d(b^*) = -ca \ko
d(c^*) = -ab \: .
\]
Then the endomorphism algebra of the cluster-tilting object $P_1\oplus P_3\oplus\tau^{-2}P_1$ in the cluster category associated with the $\R$-species
\[
\begin{tikzcd}
	\R\arrow[r,"\R"]&\R\arrow[r,"\C"]&\C
\end{tikzcd}
\]
is isomorphic to $H^0(\Pi)$.
\end{example}

\begin{remark}
All constructions and results generalise from dimension $3$ to higher dimension.
\end{remark}

\subsection*{Acknowledgement} The first-named author thanks 
Claire Amiot and Dong Yang for stimulating discussions and help with the 
references. The second-named author is supported by the
Chinese Scholarship Council (CSC, Grant No.~202006210272).
Both authors are grateful to the referee for helpful comments and
suggestions and to Damien Calaque for pointing out an error
in a previous version of the expansion \cite{KellerLiu23}
of this note.


\def\cprime{$'$} \def\cprime{$'$}
\providecommand{\bysame}{\leavevmode\hbox to3em{\hrulefill}\thinspace}
\providecommand{\MR}{\relax\ifhmode\unskip\space\fi MR }
\providecommand{\MRhref}[2]{%
  \href{http://www.ams.org/mathscinet-getitem?mr=#1}{#2}
}
\providecommand{\href}[2]{#2}

\end{document}